\documentclass[12pt]{article}

\usepackage{epic}
\usepackage{amsmath}[1996/11/01]
\usepackage{amssymb,amsthm,amsfonts,latexsym,epsfig,graphics}
 \def\beql#1#2\eeql{\begin{equation}\label{#1}#2\end{equation}}

	 \usepackage{longtable}
\usepackage{xparse} 

\ExplSyntaxOn
\NewDocumentEnvironment{flexitab}{m}
 {
  \keys_set:nn { flexitab }
   {
    type = tabular,
    stretch = 1,
    align = c,
    width = -1000pt,
    #1
   }
  \dim_compare:nTF { \l_flexitab_width_dim > -1pt }
   {
    \tl_set:Nn \l_flexitab_type_tl { tabular* }
    \tl_set:Nn \l_flexitab_start_tl { \begin{tabular*}{\l_flexitab_width_dim} }
   }
   {
    \tl_set:Nn \l_flexitab_start_tl { \begin{\l_flexitab_type_tl} }
   }
  \tl_put_right:Nx \l_flexitab_start_tl
   {
    [\l_flexitab_align_tl]{\exp_not:V \l_flexitab_columns_tl}
   }
  \tl_set_eq:NN \arraystretch \l_flexitab_stretch_tl
  \tl_use:N \l_flexitab_start_tl
 }
 {
  \str_if_eq:eeTF {\l_flexitab_type_tl}{tabular*}
   { \end{tabular*} }
   { \end{\l_flexitab_type_tl} }
 }

\keys_define:nn { flexitab }
 {
  type    .tl_set:N  = \l_flexitab_type_tl,
  align   .tl_set:N  = \l_flexitab_align_tl,
  width   .dim_set:N = \l_flexitab_width_dim,
  columns .tl_set:N  = \l_flexitab_columns_tl,
  stretch .tl_set:N  = \l_flexitab_stretch_tl,
 }
\tl_new:N \l_flexitab_start_tl
\ExplSyntaxOff

\advance \textheight by -1 \topmargin
\advance \textheight by -1 \headheight
\advance \textheight by -1 \headsep
\advance \textheight by -2 \footskip
\vsize = \textheight
\hsize = \textwidth
\DeclareMathOperator{\Comp}{Comp}

\DeclareMathOperator{\cont}{cont}
\DeclareMathOperator{\sign}{sign}

\DeclareMathOperator{\Char}{char}

\DeclareMathOperator{\Sym}{Sym}

\DeclareMathOperator{\diag}{diag}

\DeclareMathOperator{\id}{id}

\theoremstyle{plain}
\newtheorem{theorem}{Theorem}
\newtheorem{lemma}[theorem]{Lemma}

\newtheorem{proposition}[theorem]{Proposition}

\theoremstyle{remark}

\newtheorem{remark}[theorem]{Remark}

\numberwithin{theorem}{section}

\newcommand{\disj}{\stackrel{.}{\cup}}

\newcommand{\Z}{{\mathbb{Z}}}
\newcommand{\Q}{{\mathbb{Q}}}

\newcommand{\N}{{\mathbb{N}}}

\renewcommand{\em}{\sf}

\title{Symmetrizations of quadratic and hermitian forms}
\author{Gabriele Nebe\footnote{nebe@math.rwth-aachen.de}}  
\date{Lehrstuhl f\"ur Algebra und Zahlentheorie, 
RWTH Aachen University, Germany}
 
\begin{document}
\maketitle

 {\sc Abstract.} 
  \\
  The paper develops 
elementary linear algebra methods to compute the 
  determinants of the tensor symmetrizations of quadratic
  and hermitian forms over fields of good characteristic. 
  Explicit results are given for the partitions 
  $(n)$, $(1^n)$, $(2,1^{n-2})$ and $(3,1^{n-3})$ 
  as well as for all partitions of $n\leq 7$. 
  For orthogonal groups these symmetrizations are 
  not irreducible and we continue to find the 
  determinants of their irreducible constituents, 
  the refined symmetrizations, over fields of 
  characteristic 0. 
     \\
  {\sc keywords:}  symmetric group, Schur-Weyl duality, symmetrizations, 
  quadratic and hermitian forms

\section{Introduction} 

 Let $K$ be a field  and $V$ be a vector space over 
 $K$ of finite dimension, say,  $N$. 
 For $n\in \N $ the symmetric group $S_n$ acts linearly on 
 the $n$th tensor power
 $\otimes^n V$ by permuting the tensor factors thus turning $\otimes ^nV$ 
 into a $KS_n$-module.  
 To avoid trivialities we assume that $n\geq 2$ and to avoid problems 
 we assume that the characteristic of $K$ is either $0$ or $>n$. 
 It is well known that the simple $KS_n$-modules are labelled by 
 the partitions $\lambda $ of $n$. A formula for a primitive 
 idempotent $\tilde{e}_{\lambda } \in KS_n$ is given in \cite[Theorem 3.1.10]{JamesKerber}. 
 The direct summand 
 $$ \Sym_{\lambda }(V) := \tilde{e}_{\lambda } (\otimes ^n V) $$ 
 is known as the {\em $\lambda $-symmetrization of $V$}.
 Its dimension $d(\lambda , N) := \dim(\Sym_{\lambda }(V) )$ 
 is the number of semi-standard $\lambda $-tableaux with content 
 contained in $ \{ 1,\ldots , N \}$ (see \cite[Theorem 5.2.14]{JamesKerber}).

 Now let $B : V\times V \to K$ be a symmetric bilinear or 
 Hermitian sesquilinear form on $V$. 
 Then $B$ defines a respective form $\otimes ^nB$ on $\otimes ^n V $ by 
 $$\otimes ^n B(v_1\otimes \ldots \otimes v_n, w_1\otimes \ldots \otimes w_n) := 
 \prod_{i=1}^n B(v_i,w_i ) .$$
 We put $\Sym_{\lambda }(B) $ to denote the restriction of 
 $\otimes ^n B$ to $\Sym _{\lambda }(V)$. 

 One important invariant of $B$ is its {\em determinant} 
 $\det(B) $ which is defined to be the class of the 
 determinant of a Gram matrix of $B$ in 
 $K/(K^{\times })^2$ for symmetric bilinear forms and in 
 $F/N_{K/F}(K^{\times })$ for $K/F$ Hermitian forms. 

 The main result of this note is a formula for the determinant of 
 $\Sym _{\lambda }(B) $. 

 \begin{theorem}\label{mainunitary}
 $\det (\Sym _{\lambda }(B)) = c(\lambda,N) \det(B) ^{d(\lambda ,N)n/N}$,
	 where $c(\lambda , N)$ can be computed using combinatorial 
	 algorithms in $S_n$.
 \end{theorem}

 Explicit formulas for $c(\lambda ,N)$ 
 have been obtained for symmetric and exterior powers 
 of symmetric 
 bilinear forms 
 (where $\lambda = (n)$ respectively $\lambda = 1^n$)
 in \cite{symmetric} and \cite{exterior}. 
 In this paper we additionally derive 
 such formulas for the two hooks 
 $\lambda = (2,1^{n-2})$ and $\lambda = (3,1^{n-3})$. 
 For small $n$ more results are obtained by computer 
 (see Table \ref{AUS} for $n\leq 7$). 
 These are helpful to compute determinants of even degree 
 unitary characters as illustrated in Section \ref{Suz}, 
 where the symmetrizations of the $12$-dimensional unitary 
 representation of $6.Suz$ are used to obtain most of the
 determinants of the faithful, simple $\Q[\sqrt{-3}] (6.Suz)$-modules.

 For symmetric bilinear forms $B$, the symmetrizations are 
 usually not irreducible modules for the orthogonal group $O(B)$. 
 The last section investigates 
certain  $O(B)$-invariant submodules of $\Sym_{\lambda }(B)$, 
 the refined symmetrizations. For $n\leq 6$ a
  table that can be used to compute 
  their determinants is given in Section \ref{AUSorth}.
  As an application we obtain some orthogonal determinants 
  for the sporadic simple Conway group $Co_1$, that were not
  contained in the database described in \cite[BBBNP23]{survey} yet.

 I thank Thomas Breuer for motivation, helpful comments and 
 pointing out the reference \cite{Frame}. 

 \section{Notation}

 We use the standard notation as given in the textbook \cite{JamesKerber}. 
 A {\em partition} $\lambda = ( \lambda_1,\ldots , \lambda _k)$ of 
a natural number $n$ is a sequence of integers $\lambda _i$ with 
 $\lambda _1\geq \lambda _2 \geq \ldots \geq \lambda _k$ 
 such that $\sum_{i=1}^k \lambda_i = n$.
 It is visualised by its {\em Young diagram} with 
 $k$ rows and $\lambda _i $ boxes in row $i$. 
A {\em Young tableau} $t^{\lambda }$ is obtained by 
labelling these $n$ boxes of the Young diagram by the $n$ numbers 
in $\{1,\ldots , n\}$. 
Such a Young tableau defines two  partitions of $\{ 1,\ldots ,n \}$, 
the horizontal partition $\disj _{j=1}^k P_j = \{ 1,\ldots , n\}$ 
where $P_j$ consists of the $\lambda _j$ labels in row number $j$ 
and a corresponding vertical partition given by the 
labels in the $\lambda _1$ columns of the Young tableau.
The horizontal group $H^{\lambda }$ of $t^{\lambda }$ is the 
stabiliser in $S_n$ of all sets in the horizontal partition and 
the vertical group $G^{\lambda }$ the stabiliser of all sets in the 
vertical partition.  

With this notation the formula for a primitive idempotent 
$\tilde{e}_{\lambda } $ 
from \cite[Theorem 3.1.10]{JamesKerber} is 
$\tilde{e}_{\lambda }= \frac{d_{\lambda }}{n!} e_{\lambda} $ with
$$e_{\lambda } = 
\sum _{\sigma \in G^{\lambda }} \sum _{\rho \in H^{\lambda}}  
\sign(\sigma ) \sigma \rho .$$

An {\em $N$-semi-standard Young tableau} $t$ of shape $\lambda $ 
is given by filling elements in $\{1,\ldots , N \}$ into the boxes
of the Young diagram of $\lambda $ such that they are non-decreasing along the 
rows and strictly increasing along the columns. 
The {\em content} $\cont(t)$ of $t$ is the 
multi-set of its entries in $\{ 1,\ldots ,N \}$. 
Put 
$$T(\lambda ,N ):= \{ t \mid t \mbox{ in an $N$-semi-standard Young tableau
of shape $\lambda $ } \} .$$

The following result is well known: 

\begin{remark}
Let  $d(\lambda , N) := |T(\lambda ,N) |$ 
denote the number of $N$-semi-standard Young tableaux
of shape $\lambda $.
Then each element of $\{ 1,\ldots , N \}$ 
occurs with the same multiplicity in the union of 
the contents of the elements of $T(\lambda ,N) $ and hence 
$$\bigcup _{t \in T(\lambda , N)} \cont (t) 
= \{ 1^{d(\lambda ,N) n/N} , \ldots , N^{d(\lambda ,N) n/N } \} .$$
\end{remark}

 \section{Proof of the main theorem} 

 Let $(v_1,\ldots , v_N)$ be an orthogonal basis of $(V,B)$ and 
 put $a_i:=B(v_i,v_i)$. 
 Then $\det (B) $ is represented by $a_1\cdots a_n$. 

 \begin{remark} 
	 The pure tensors $v_{i_1} \otimes \ldots \otimes 
	 v_{i_n} $ form an orthogonal basis of 
	 $(\otimes ^n V , \otimes ^n B)$ with 
	 $$\otimes ^n B ( v_{i_1} \otimes \ldots \otimes v_{i_n} , v_{i_1} \otimes \ldots \otimes v_{i_n} ) = \prod _{j=1}^n a_{i_j} .$$
 \end{remark}

 Fix a partition $\lambda $ of $n$ and fix a Young tableau 
 $t^{\lambda }$ to label the positions in $\otimes ^n V$. 
Then any $t\in T(\lambda , N)$ defines a 
 pure tensor  
 $$v(t) := w_1\otimes \ldots \otimes w_n \in \otimes ^n V$$
 where $w_i = v_j$, if the box that is labelled by 
 $i$ in $t^{\lambda }$ has content $j$ in $t$.

 \begin{remark} 
	 Let $t\in T(\lambda, N) $ and $C:=\cont(t)$ be its content.
	 Then
	 $$\otimes ^n B(v(t) , v(t) )  = 
	 \prod _{i\in C } a_i =: q(C)$$ 
	 only depends on $C$.
 \end{remark}

 \begin{lemma} (\cite[Section 4.3]{James}, \cite[Section 5.2]{JamesKerber}) \label{basis} 
	 A basis of $\Sym _{\lambda }(V)$ is given by 
	 $$\{ e_{\lambda } (v(t) ) \mid 
	 t \in T(\lambda , N ) \} .$$
 \end{lemma}

 As the summands in $e_{\lambda } (v(t) ) $ are 
 (up to a sign) pure tensors of the basis elements $v_i$ with the 
 same multi-set of indices, the value  
 $$\otimes ^n B(e_{\lambda } (v(t) ) , e_{\lambda } (v(t) ) ) = c(t) q(C) $$
 is a constant multiple of $q(C)$ where $C=\cont(t)$.

 The following trivial remark is fundamental for the computations: 
 \begin{remark}\label{oooo} 
For $t,s\in T(\lambda,N)$ we get
 $$\otimes ^n B (e_{\lambda } (v(t) ) , e_{\lambda } (v(s) ) )  
 = 0 $$ 
 unless $\cont (t) = \cont (s)$.  
 \end{remark} 

 Ordering the basis from Lemma \ref{basis} according to the 
 content $C$ of $t$ we hence obtain a Gram matrix of 
 $\Sym_{\lambda }(B)$ as a block diagonal matrix 
 $$\diag ( q(C) X_C : C \in \cont(T(\lambda, N) ) ) $$ 
 and hence Theorem \ref{mainunitary} follows with 
 $$c(\lambda ,N) = \prod _{C\in \cont(T(\lambda ,N))}  \det (X_C) .$$

For computing the Gram matrices $X_C$ I wrote a small ad hoc program.
Up to $n=7$, the results are given in Table \ref{AUS}.

\section{Small examples}

\subsection{The partition $(2,1)$} 

This partition is the smallest case that it not yet available in the 
literature. 
To illustrate the algorithm used to compile Table \ref{AUS} 
I explain these computations in detail. 
The formula in \cite[Theorem 5.2.14]{JamesKerber} yields 
$$\dim (\Sym_{(2,1)} (V)) = \frac{1}{3} N(N-1)(N+1) , 
\mbox{ where } N= \dim(V) .$$
So 
$$\det(\Sym_{(2,1)}(B) )= \det(B) ^{(N-1)(N+1)} c((2,1),N) $$
for some $c((2,1),N)$. 

\begin{proposition}
	$c((2,1),N) = 3^{{N}\choose{3}} z^2 $
	for some $z\in \N $ with prime divisors in $\{ 2,3 \}$.
\end{proposition}

\begin{proof}
	We give the 
	$N$-semi-standard Young tableaux of shape $(2,1)$ 
	by listing the 
	entries in positions $((1,1), (1,2), (2,1) ) $ in this ordering.

Choosing $a<b<c $, $a,b,c\in \{ 1,\ldots , N \}$ there are 
	4 sorts of $N$-semi-standard $(2,1)$-tableaux: 
	$$\{ (a,a,b), (a,b,b), (a,b,c) ,(a,c,b) \} $$
	Then 
	$$\begin{array}{lcl}
		{e}_{(2,1)} (a,a,b) &=& 2((a,a,b)-(b,a,a))  \\
		{e}_{(2,1)} (a,b,b) &=& (a,b,b)-(b,b,a)  \\
		{e}_{(2,1)} (a,b,c) &=& (a,b,c)-(c,b,a) + (b,a,c) - (c,a,b)  \\
		{e}_{(2,1)} (a,c,b) &=& (a,c,b)-(b,c,a)+(c,a,b)-(b,a,c)  \\
	\end{array}
	$$
	giving rise to 
	$$
	\begin{flexitab}{
			type = array,
			stretch=1.3,
			columns=lcl
			} 
		X_{a^2b} = (8) & \mbox{ with multiplicity } & {{N}\choose{2}} \\ 
		X_{ab^2} = (2) & \mbox{ with multiplicity } & {{N}\choose{2}} \\ 
		X_{abc} =  \left( \begin{array}{rr} 4 & -2 \\ -2 & 4 \end{array}\right)   & \mbox{ with multiplicity } & {{N}\choose{3}} 
	\end{flexitab} 
	$$
	yielding
	 $$\det(\Sym_{(2,1)} (B)) = 16^{{N}\choose{2}} 12^{{N}\choose{3}} 
	 \det(B)^{(N-1)(N+1)} .$$
\end{proof}

\subsection{Two partitions of 4} \label{part4}

\begin{proposition}
	$$\dim(\Sym_{(3,1)} (V)) = \frac{1}{8} (N+2)(N+1)N(N-1) = 
	3 {{N+2}\choose{4}}  $$ 
	and, up to squares, 
	$$\det(\Sym_{(3,1)}(B)) =  2^{{N}\choose{3}} \det(B)^{(N+2)(N+1)(N-1)/2}.$$
\end{proposition}

\begin{proof}
	For $\lambda = (3,1) $ we compute the relevant $X_C$, $C\in \cont(T(\lambda , N) $ as in the previous section and obtain: 
	$$
	\begin{flexitab}{
                        type = array,
                        stretch=1.3,
                        columns=lcl
                        }
 X_{a^3b} = (72) & \mbox{ with multiplicity } & {{N}\choose{2}} \\
X_{a^2b^2} = (16) & \mbox{ with multiplicity } & {{N}\choose{2}} \\
X_{ab^3} = (8) & \mbox{ with multiplicity } & {{N}\choose{2}} \\
	       X_{a^2bc} = \begin{pmatrix} 24 & -8 \\ -8 & 24 \end{pmatrix}  & \mbox{ with multiplicity } & {{N}\choose{3}} \\
		       X_{ab^2c} =  \begin{pmatrix} 24 & -8 \\ -8 & 8 \end{pmatrix}  & \mbox{ with multiplicity } & {{N}\choose{3}} \\
			       X_{abc^2} =  \begin{pmatrix} 24 & -8 \\ -8 & 8 \end{pmatrix} & \mbox{ with multiplicity } & {{N}\choose{3}} \\
				       X_{abcd} =  \begin{pmatrix} 12 & -4 & -4  \\ -4 & 12 & -4 \\ -4 & -4 & 12 \end{pmatrix}  & \mbox{ with multiplicity } & {{N}\choose{4}} 
        \end{flexitab}
        $$
	Up to squares $c(\lambda,N) = \det(X_{a^2bc}) ^{{N}\choose{3}} $.
\end{proof}

\begin{proposition}
	$$\dim(\Sym_{(2,2)} (V)) = \frac{1}{12} (N^4-N^2) $$ 
	and, up to squares, 
	$$\det(\Sym_{(2,2)} (B)) = 
	3^{{N}\choose{4}} 2^{{N}\choose{3}} \det(B)^{N(N+1)(N-1)/3}$$
\end{proposition}

\begin{proof}
	For $\lambda = (2,2) $ we obtain 
	$$
	\begin{flexitab}{
                        type = array,
                        stretch=1.3,
                        columns=lcl
                        }
 X_{a^2b^2} = (64) & \mbox{ with multiplicity } & {{N}\choose{2}} \\
 \\
 X_{a^2bc} = (32) & \mbox{ with multiplicity } & {{N}\choose{3}} \\
 \\
X_{ab^2c} = (8) & \mbox{ with multiplicity } & {{N}\choose{3}} \\
\\
X_{abc^2} = (32) & \mbox{ with multiplicity } & {{N}\choose{3}} \\
\\
	       X_{abcd} = \begin{pmatrix} 16 & -8  \\ -8 & 16 \end{pmatrix}  & \mbox{ with multiplicity } & {{N}\choose{4}} 
        \end{flexitab}
        $$
	Up to squares $c(\lambda ,N) = \det(X_{ab^2c})^{{N}\choose{3}} 
	\det(X_{abcd}) ^{{N}\choose{4}}  $.
\end{proof}

\subsection{The partitions $(n)$ and $(1^n)$} 

The module $\Sym_{(1^n)}(V) = \Lambda ^n(V)$ is better known 
as the $n$-th exterior power of $V$, whereas 
$\Sym_{(n)}(V) $ is the $n$-th symmetric power of $V$.

\begin{proposition} (see also \cite{exterior}) 
$\dim(\Sym_{(1^n)}(V) ) =  {{N}\choose{n}} $ 
	and $$\det(\Sym_{(1^n)}(B)) = 
	  (n!)^{{{N}\choose{n}} } \det(B) ^{{{N-1}\choose{n-1}}} .$$
\end{proposition}

\begin{proof}
	Let $\lambda := (1^n)$. 
	Then the horizontal group $H^{\lambda } = \{1 \}$ 
	and the vertical group $G^{\lambda } = S_n$.
	Moreover the $N$-semi-standard tableaux of shape 
	$\lambda $ are in bijection with the $n$-element subsets of $\{ 1,\ldots, N\}$, in particular 
	$$\dim (\Sym _{\lambda }(V)) = 
	|T(\lambda , N) |= {{N}\choose{n}} $$
	and any element $t\in T(\lambda ,N)$ is uniquely determined by its 
	content $C=\cont(t)$. So $X_C$ is a scalar and equals to the 
	length of the $G^{\lambda }$-orbit of $t$, which shows that 
	$$\otimes ^nB (e_{\lambda }(v(t)) , e_{\lambda }(v(t))) = 
	 n! q(C) .$$
	 So 
	 $$\det(\Sym_{(1^n)}(B)) = (n!)^{{{N}\choose{n}} } \prod _{1\leq i_1 < \ldots < i_n \leq N } a_{i_1}\cdots  a_{i_n} .$$
	 As each element $i\in \{1,\ldots , N\}$ occurs 
 ${N-1}\choose{n-1}$ times in an $n$-element subset of $\{1,\ldots , N\}$  
	the latter product is just 
	$\prod_{i=1}^N a_i^{{{N-1}\choose{n-1}}} = \det(B) ^{{{N-1}\choose{n-1}}} .$
\end{proof}

	To state the result for  $\lambda = (n)$ we introduce the set 
	$$\Comp(n,k) := \{ (x_1,\ldots , x_k) \in \Z_{>0}^k \mid 
	x_1 + \ldots + x_k = n \} $$ 
	 of all compositions of $n$. 

\begin{proposition} (see also \cite[Proposition 3.9]{symmetric})
$\dim(\Sym_{(n)}(V) ) =  {{N+n-1}\choose{n}} $ 
	and $$\det(\Sym_{(n)}(B)) = 
	\prod_{k=1}^n \prod _{(x_1,\ldots , x_k) \in \Comp(n,k) }
           (\frac{n!}{(x_1!)\cdots (x_k!) })^{{N}\choose{k} }
	   \det(B) ^{{N+n-1}\choose{n-1}}  
	   .$$
\end{proposition}

\begin{proof}
	For $\lambda = (n)$ 
	 the vertical group of $\lambda $ is trivial and 
 the horizontal group $H^{\lambda } = S_n$. 
	Now the $N$-semi-standard tableaux of shape 
	$\lambda $ are in bijection with the $n$-element multi-subsets of $\{ 1,\ldots, N\}$, in particular 
	$$\dim (\Sym _{\lambda }(V)) = 
	|T(\lambda , N) |= {{N+n-1}\choose{n}} .$$
	Again any 
	element $t\in T(\lambda ,N)$ is uniquely determined by its 
	content $C = \cont(t) $ and $X_C$ is a scalar. 
	Instead of working with ${e}_{\lambda } v(t)$ 
	we divide this vector by the order of the stabiliser of $t$ 
	and hence just work with the orbit sum 
	$$\sum  \{ v(\sigma(t)) \mid \sigma \in S_n \} .$$
	This results in multiplying $X_C$ be an integral  square. 
	If $C= \{ n_1^{x_1},\ldots , n_k ^{x_k } \}$  with 
	$n_1\leq \ldots , \leq n_k$ then 
	$X_C = \frac{n!}{(x_1!)\cdots (x_k!)} $.
	Therefore
	$$c(\lambda ,N) = \prod_{k=1}^n \prod_{(x_1,\ldots , x_k) \in \Comp(n,k) }
	   (\frac{n!}{(x_1!)\cdots (x_k!) })^{{N}\choose{k} }.$$
\end{proof}

\subsection{A result for hooks} 

Another situation, where $X_C$ can be computed for general $\lambda $ and $N$,
is if $C$ is a subset of $\{ 1,\ldots , N\} $ or cardinality $n$, i.e. every 
element of the multi-set $C$ occurs in $C$ with multiplicity $1$.

\begin{lemma}\label{hook} 
	Assume that $\lambda = (\ell , 1^{n-\ell})$ is a hook.
	If $C$ is  multiplicity free then 
	$\langle e_{\lambda } v(t)  \mid t \in T(\lambda ,N) , \cont(t) = C \rangle \cong S^{\lambda '}$ is the irreducible representation of 
	$S_n$ isomorphic to the Specht module 
	associated with the transposed partition $\lambda '=(n-\ell+1,1^{\ell-1})$
	and $X_C$ is a Gram matrix of an $S_n$-invariant symmetric 
	bilinear form
	on $S^{\lambda '}$ with
$$\det(X_C) = ((\ell -1) ! (n-\ell +1) ! )^{{{n-\ell}\choose{\ell-1}}}  n^{{{n-1}\choose{\ell-1}}} .$$
\end{lemma}

\begin{proof}
	Let $C\subseteq \{ 1,\ldots , N\}$ with $|C|=n$ and let $a$ be 
	the minimal element of $C$. 
	Then the semi-standard $\lambda $-tableaux with content $C$ are
	the ${{n-\ell}\choose{\ell-1}} $ elements of 
	$$\{ t_S \mid S \subseteq C\setminus \{ a \} , |S| = \ell -1 \} $$
	where the elements of the first row of $t_I$ 
	form the set $S\cup \{ a\} $. 
	As all summands in $e_{\lambda } v(t_S) $ are distinct we have 
	$\otimes ^n B (e_{\lambda } v(t_S) , e_{\lambda } v(t_S) ) = 
	\ell ! (n-\ell +1) ! .$ 
	For $R\neq S\subseteq C\setminus \{ a \}  $ the 
	common summands of $e_{\lambda } v(t_R) $ and $e_{\lambda } v(t_S) $ 
	are those $v(t)$, where the set of elements in the 
	first column of $t$ are identical, i.e. the last $\ell -1$ elements
	of the first row are $\{a \} \cup (R\cap S)$, i.e. we
	only obtain non-zero inner product if $|R\cap S| = \ell -2$. 
	So 
	$$\otimes ^n B (e_{\lambda } v(t_S) , e_{\lambda } v(t_R) ) = 
	\left\{ \begin{array}{ll} 	
		\ell  ! (n-\ell +1) ! & \mbox{ if } S=R  \\
		\pm (\ell -1) ! (n-\ell +1) ! & \mbox{ if } |R\cap S| = \ell -2  \\
		0 & \mbox{ if } |R\cap S| < \ell -2  
	\end{array} \right) $$
	where the sign depends on the sign of the permutation $\pi_{R,S}$
	mapping $[x,b,c,d,\ldots ]$ to $[y,b,c,d,\ldots ] $ where
	$S\setminus R = \{x \} $ and $R\setminus S = \{ y \} $.

	To compute the $(\ell -1)$st exterior power of the root lattice
	$A_{n-1}$ we choose a basis 
	$$A_{n-1} = \langle b_2 := e_1-e_2,\ldots , e_1-e_n =: b_n \rangle $$
	such that the Gram matrix of $(b_2,\ldots ,b_n)$ is 
	$I_{n-1} + J_{n-1}$ of determinant $n$. 
	To an $\ell -1$ element subset $R=\{i_1,\ldots , i_{\ell-1}\}$ of $C \setminus \{ 1\} = \{2,\ldots , n\}$
	we associate the basis vector 
	$b_R = b_{i_1}\wedge \ldots \wedge b_{i_{\ell -1}} $. 
	So in $\Lambda ^{\ell -1} (A_{n-1})$ we obtain the inner product of 
	the subsets $R$ and $S$ as the $\ell-1 \times \ell -1 $-minor
	$$\det (I_{n-1} + J_{n-1}) _{R\times S} = 
	\pm \det (\left( \begin{array}{cc} 0 & 0 \\ 0 & I_{R\cap S}  \end{array} 
		\right) + J_{\ell-1} ) 
		=\left\{ \begin{array}{ll} 
			\ell & \mbox{ if } R=S \\
			\sign(\pi_{R,S})  & \mbox{ if } |R\cap S| = \ell-2 \\
			0 & \mbox{ if } |R\cap S| < \ell-2 .
		\end{array} \right. $$
		So the lattice
		$$(\Z^{{{n-\ell}\choose{\ell-1}}} , ((\ell -1) ! (n-\ell +1) ! )^{-1} X_C ) \cong \Lambda ^{\ell-1} (A_{n-1}) $$ 
		and hence 
		$$\det(X_C) = ((\ell -1) ! (n-\ell +1) ! )^{{{n-\ell}\choose{\ell-1}}}  n^{{{n-1}\choose{\ell-1}}} .$$
\end{proof}

\subsection{The partitions $(2,1^{n-2})$} 

\begin{proposition}
	Let $\lambda := (2,1^{n-2})$. 
	Then $$\dim(\Sym_{\lambda}(V)) = (n-1) {{N}\choose{n}} +
	(n-1) {{N}\choose{n-1}} 
	=(n-1) {{N+1}\choose{n}} $$ and 
	$$c(\lambda ,N) = n^{{N}\choose{n}} ((n-1)!)^{(n-1){{N+1}\choose{n}}} 
	.$$
\end{proposition}

\begin{proof}
	For $t\in T(\lambda ,N)$ the content $C$ of $t$ has either 
	$n-1$ or $n$ distinct elements. 
	If one element occurs twice in $C$, then this is either the 
	minimal element of $C$ and 
	 $t=(a,a,b,c,d,\ldots )$ such that 
	$$e_{\lambda } v(t) = 2 \sum_{g\in S_{n-1}} \sign(g) v(gt) $$
	or it is not 
	the minimal element and the tableau is 
	$s=(a,x,b,c,d,\ldots )$ where $x\in \{ b,c,d,\ldots \}$ and then 
	$$e_{\lambda } v(s) = \sum_{g\in S_{n-1}} \sign(g) v(gs) .$$ 
	In both cases the tableau is uniquely determined by its content. 
	We compute 
	$$B(e_{\lambda }v(s), e_{\lambda } v(s) ) = 
	B(1/2 e_{\lambda }v(t), 1/2 e_{\lambda } v(t) ) =  (n-1)! q(C) .$$
	This gives a contribution of $(n-1) {{N}\choose{n-1}}$ to 
	$\dim(\Sym_{\lambda}(V))$ and  of
	$$(n-1)!^{(n-1) {{N}\choose{n-1}}}  \mbox{ to } 
c(\lambda ,N).$$ 
	If every element in $C=\{ a,b,c,d ,\ldots \} $ 
	occurs with multiplicity 1, 
	then Lemma \ref{hook} shows that $X_C\cong  (n-1)! A_{n-1} $.
	In total these $C$ contribute  $(n-1){{N}\choose{n}} $ 
	to the dimension and 
	$$(n (n-1)!^{(n-1)} ) ^{{N}\choose{n}} \mbox{ to } 
c(\lambda , N).$$
\end{proof}

\subsection{The partitions $(3,1^{n-3})$} 

\begin{proposition}
	Let $\lambda := (3,1^{n-3})$. 
	Then $$\dim(\Sym_{\lambda}(V)) = {{n-1}\choose{2}} 
	( {{N}\choose{n}} + 2 {{N}\choose{n-1}} + {{N}\choose{n-2}})
	$$ and 
	$$c(\lambda ,N) = (n-2)!^x  2^y n^z $$ 
	where 
	$$ 
		\begin{flexitab}{
                        type = array,
                        stretch=1.3,
                        columns=lcl
                        }
		x & = & 
		({{N}\choose{n-2}}+{{N}\choose{n}} )  {{n-1}\choose{2}} 
		(n-2) {{N}\choose{n-1}}   \\
		y & = & 
		{{N}\choose{n-2}} {{n-2}\choose{2}}  +
		(n-3) {{N}\choose{n-1}}   +
		{{N}\choose{n}} )  {{n-1}\choose{2}}  
		\\
		z & = & {{N}\choose{n-1}} + (n-2) {{N}\choose{n}}  
	. \end{flexitab} $$
\end{proposition}

\begin{proof}
	For $t\in T(\lambda ,N)$ the content $C$ of $t$ has either 
	$n-2$, $n-1$ or $n$ distinct elements. 

	If one element occurs three times in $C$, then this is either the 
	minimal element of $C$ and 
	 $t=(a,a,a,b,c,d,\ldots )$ such that 
	$$e_{\lambda } v(t) = 6 \sum_{g\in S_{n-2}} \sign(g) v(gt) $$
	or it is not 
	the minimal element and the tableau is 
	$s=(a,x,x,b,c,d,\ldots )$ where $x\in \{ b,c,d,\ldots \}$ and then 
	$$e_{\lambda } v(s) = 2\sum_{g\in S_{n-2}} \sign(g) v(gs)  .$$
	In both cases the tableau is uniquely determined by its content. 
	We compute 
	$$B(1/2 e_{\lambda }v(s), 1/2 e_{\lambda } v(s) ) = 
	B(1/6 e_{\lambda }v(t), 1/6 e_{\lambda } v(t) )  =
	(n-2)! q(C) 
	$$
So we obtain a contribution of ${{N}\choose{n-2}} (n-2)$ to the dimension and of
	$$(n-2)!^{{{N}\choose{n-2}} (n-2) } \mbox{ to } c(\lambda , N).$$
	If two elements occur twice in $C$, then either one of them is the
        minimal element of $C$ and
         $t=(a,a,x,b,c,d,\ldots )$ 
         where $x\in \{ b,c,d,\ldots \}$. Then
	$$e_{\lambda } v(t) = 2(\sum_{g\in S_{n-2}} \sign(g) v(gt)  
	+
	\sum_{g\in S_{n-2}} \sign(g) v(gt')) 	$$
	where $t'= (a,x,a,b,c,d,\ldots )$.
	Or $s=(a,x,y,b,c,d,\ldots )$ where 
	$x<y\in C \setminus \{ a\}$
	and 
	$$e_{\lambda } v(s) = \sum_{g\in S_{n-2}} \sign(g) v(gs)  
	+
\sum_{g\in S_{n-2}} \sign(g) v(gs') 	$$
	where $s'= (a,y,x,b,c,d,\ldots )$.
        We compute
        $$B(1/2 e_{\lambda }v(t), 1/2 e_{\lambda } v(t) ) =
        B( e_{\lambda }v(s),  e_{\lambda } v(s) ) =
        2 (n-2)! q(C)
        $$
So we obtain a contribution of ${{N}\choose{n-2}} {{n-2}\choose{2}} $ to the 
dimension and  of 
	$$(2(n-2)!)^{{{N}\choose{n-2}} {{n-2}\choose{2}} } 
	\mbox{ to } c(\lambda ,N).$$
	If $C$ consists of $n-1$ distinct elements, then one 
	element of $C$ occurs twice. 
	If this is the minimal element of $C$, then 
	 $t=(a,a,x,b,c,d,\ldots )$ 
         where $x\not\in \{ a,b,c,d,\ldots \}$. Then
        $$e_{\lambda } v(t) = 2(\sum_{g\in S_{n-2}} \sign(g) v(gt)  
        + \sum_{g\in S_{n-2}} \sign(g) v(gt')
        + \sum_{g\in S_{n-2}} \sign(g) v(gt'')
	)    $$
	where $t'= (a,x,a,b,c,d,\ldots )$ and $t'' = (x,a,a,b,c,d,\ldots )$.
	We compute
	$$B(1/2 e_{\lambda }v(t), 1/2 e_{\lambda } v(t) ) = 3 (n-2) ! .$$
	There are in total $n-2 = |\{ x,b,c,d,\ldots \} |$ 
	semi-standard $\lambda $ tableaux with the same content. 
	The inner product of distinct such tableaux is $\pm (n-2) ! $
	depending on the product of the signs of the permutations sorting the 
	respective $[x,b,c,d,\ldots ] $.
	So the determinant of $X_C$ is 
	$$(n-2)!^{n-2} (\det (2I_{n-2}+J_{n-2})) = 
	(n-2)!^{n-2} 2^{n-3} n $$
	adding $(n-2) {{N}\choose{n-1}} $ to the dimension and 
	a factor 
	$$((n-2)!^{n-2} 2^{n-3} n)^{{N}\choose{n-1}} \mbox{ to the determinant.} $$
	If a non-minimal element occurs twice in $C$, then the tableau 
	is either 
	$$t=(a,x,x,b,c,d,\ldots )\mbox{ or } 
	s=(a,\min(x,y),\max(x,y),b,\ldots , x , \ldots ) $$
	We get 
	$$B(1/2 e_{\lambda }v(t), 1/2 e_{\lambda } v(t) ) = 3 (n-2) ! ,\ 
	B( e_{\lambda }v(s),  e_{\lambda } v(s) ) = 4 (n-2) ! ,$$  
	In total there are $n-2$ such semi-standard 
	tableaux with the same content $C$,
	one tableau $t$ and 
	$(n-3)$ possibilities for tableaux of type $s$
	depending on the choice of $y\in \{b,c,d,\ldots \}$. 
	For such $s,s'$ we compute
	$$B(1/2 e_{\lambda }v(t),  e_{\lambda } v(s) ) = \pm 2 (n-2) ! = 
	B( e_{\lambda }v(s),  e_{\lambda } v(s') ) .$$  
	where the sign is the signum of the permutation 
	mapping $[x,b,c,d,\ldots ]$ to $[y,b,c,d, \ldots ]$ 
	respectively $[y,b,c,d,\ldots ]$ to $[y',b,c,d,\ldots ]$. 
	So after multiplying those basis vectors 
	$e_{\lambda } v(s) $ by $-1$  that have $B(1/2 e_{\lambda }v(t),  e_{\lambda } v(s) ) = - 2 (n-2) !$  we obtain 
	$$X_C = (n-2)! \left(\begin{array}{cccc} 3 & 2 & \ldots & 2  \\
		2 & 4 & \ldots & 2 \\
		\vdots & \ddots & \ddots & \vdots \\
	2 & \ldots & 2 & 4 \end{array} \right) $$
	of determinant $(n-2)!^{n-2} 2^{n-3} n .$
	Such $C$ contribute in total 
	$(n-2) {{N}\choose{n-1}}$ to the dimension and 
	$$
	((n-2)!^{n-2} 2^{n-3} n )^{{N}\choose{n-1}} \mbox{ to } c(\lambda ,N).$$
	If every element in $C=\{ a,b,c,d ,\ldots \} $ 
	occurs with multiplicity 1, 
	then Lemma \ref{hook} yields that 
	$X_C \cong 2(n-2)! \Lambda ^2 (A_{n-1})  $ 
	has determinant 
	$$(2(n-2)!)^{{{n-1}\choose{2}}} n^{n-2} .$$
	So we obtain an additive contribution of 
	${{n-1}\choose{2}} {{N}\choose{n}} $ to the dimension 
	and a multiplicative contribution of 
	$$((2(n-2)!)^{{{n-1}\choose{2}}} n^{n-2})^{{N}\choose{n}} \mbox{ to the 
	determinant.} $$
\end{proof}

\section{Determinants of symmetrizations} \label{AUS}

This section gives  tables of dimensions and determinants of 
$(\Sym_{\lambda }(V), \Sym_{\lambda }(B)) $ 
for partitions $\lambda $ of $n$ and all $n\leq 7$.
Here $(V,B)$ is either a non-degenerate symmetric bilinear or  Hermitian 
space over a field of characteristic not dividing $n!$. 
For fixed partition $\lambda $ of $n$ the dimension 
 $d(\lambda ):=\dim (\Sym_{\lambda }(V) )$ is a polynomial in $N:=\dim(V)$.
We also give $c(\lambda )$ such that 
$\det (\Sym_{\lambda }(B)) = c(\lambda ) \det(B)^{d(\lambda ) n/N }$,
up to squares.

$$
\begin{flexitab}{
                type = array,
                stretch = 1.4,
		columns = |cccc|
                }
		\hline
n & \lambda  &  d(\lambda ) & c(\lambda )  \\ 
\hline 
2 & (2) & N(N+1)/2 & 2^{{N}\choose{2}} \\
 & (1^2) & N(N-1)/2 & 2^{{N}\choose{2}} \\
\hline
3 & (3) & N(N+1)(N+2)/6 & 6^{{N}\choose{3}} \\
 & (2,1) & N(N-1)(N+1)/3 & 3^{{N}\choose{3}} \\
 & (1^3) & N(N-1)(N-2)/6 & 6^{{N}\choose{3}} \\
\hline
4 & (4) & N(N+1)(N+2)(N+3)/24 & 2^{{{N}\choose{2}} + {{N}\choose{4}}}
 3^{{{N}\choose{2}} + {{N}\choose{3}} + {{N}\choose{4}}}
 \\
 & (3,1) & N(N-1)(N+1)(N+2)/8 & 2^{{N}\choose{3}} \\
 & (2^2) & N^2(N-1)(N+1)/12 & 2^{{N}\choose{3}}3^{{{N}\choose{4}}} \\
 & (2,1^2) & N(N-1)(N-2)(N+1)/8 & 6^{{{N}\choose{3}}+{{N}\choose{4}}} \\
 & (1^4) & N(N-1)(N-2)(N-3)/24 & 6^{{{N}\choose{4}}} \\
\hline
	5 & (5) & N(N+1)(N+2)(N+3)(N+4)/120 & 6^{{{N}\choose{3}} + {{N}\choose{5}}}
	5^{{{N}\choose{5}}}
 \\
	& (4,1) & N(N-1)(N+1)(N+2)(N+3)/30 & 3^{{N}\choose{3}} 
	5^{{{N}\choose{5}}} \\
	& (3,2) & N^2(N-1)(N+1)(N+2)/24 & 3^{{N}\choose{5}} \\
	& (3,1^2) & N(N-1)(N-2)(N+1)(N+2)/20 & 2^{{{N}\choose{3}}}5^{{{N}\choose{5}}} \\
	& (2^2,1) & N^2(N-1)(N-2)(N+1)/24 & 3^{{{N}\choose{3}}}6^{{{N}\choose{5}}} \\
	& (2,1^3) & N(N-1)(N-2)(N-3)(N+1)/30 & 5^{{{N}\choose{5}}} \\
	& (1^5) & N(N-1)(N-2)(N-3)(N-4)/120 & 30^{{{N}\choose{5}}} \\
\hline
\end{flexitab}
$$

$$
\begin{flexitab}{
		type = array,
		stretch = 1.4, 
		columns = |ccc|
		}
	\hline
 \lambda  &  d(\lambda ) & c(\lambda )  \\ 
\hline 
	 (6) & N(N+1)(N+2)(N+3)(N+4)(N+5)/720 &
	2^{{{N}\choose{5}}}
	3^{{{N}\choose{3}}}
	5^{{{N}\choose{2}}+ {{N}\choose{5}} + {{N}\choose{6}}}
 \\
	 (5,1) & N(N-1)(N+1)(N+2)(N+3)(N+4)/144 & 2^{{N}\choose{3}} 
	3^{{{N}\choose{2}}+ {{N}\choose{5}} + {{N}\choose{6}}}
	 \\
	 (4,2) & N^2(N-1)(N+1)(N+2)(N+3)/80 & 
	 2^{{{N}\choose{2}} + {{N}\choose{6}}} 5^{{N}\choose{5}} \\
	 (4,1^2) & N(N-1)(N-2)(N+1)(N+2)(N+3)/72 & 6^{{{N}\choose{3}}}3^{{{N}\choose{5}}} \\
	 (3^2) & N^2(N-1)(N+1)^2(N+2)/144 & 
	 2^{{{N}\choose{2}} + {{N}\choose{6}}} 3^{{{N}\choose{3}} + {{N}\choose{5}}} \\
	 (3,2,1) & N^2(N-1)(N-2)(N+1)(N+2)/45 & 3^{{{N}\choose{3}}}5^{{{N}\choose{5}}} \\
	 (3,1^3) & N(N-1)(N-2)(N-3)(N+1)(N+2)/72 & 3^{{{N}\choose{5}}} \\
	 (2^3) & N^2(N-1)^2(N-2)(N+1)/144 & 3^{{{N}\choose{5}}} \\
	 (2^2,1^2) & N^2(N-1)(N-2)(N-3)(N+1)/80 & 3^{{{N}\choose{5}}} 30^{{{N}\choose{6}}}
	\\
	 (2,1^4) & N(N-1)(N-2)(N-3)(N-4)(N+1)/144 & 5^{{{N}\choose{6}}} 30^{{{N}\choose{5}}}
	\\
	 (1^6) & N(N-1)(N-2)(N-3)(N-4)(N-5)/720 & 5^{{{N}\choose{6}}} \\
\hline
\end{flexitab}
$$

$$
\begin{flexitab}{
                type = array,
                stretch = 1.4,
                columns = |ccc|
                }
	\hline
 \lambda  &  d(\lambda ) & c(\lambda )  \\ 
\hline 
	 (7) & {{N+6}\choose{7}}  &
	6^{{{N}\choose{5}}}
	5^{{{N}\choose{5}} + {{N}\choose{7}}}
	7^{{{N}\choose{3}}+ {{N}\choose{5}} + {{N}\choose{7}}}
 \\
	 (6,1) & \text{\scriptsize{$N(N-1)(N+1)(N+2)(N+3)(N+4)(N+5)/840$}} & 
	3^{{{N}\choose{3}}+ {{N}\choose{5}} }
	 5^{{N}\choose{3}} 
	7^{{{N}\choose{3}}+ {{N}\choose{5}} + {{N}\choose{7}}}
	 \\
	 (5,2) & \text{\scriptsize{$N^2(N-1)(N+1)(N+2)(N+3)(N+4)/360$}} & 
	 2^{{{N}\choose{3}} } 
	3^{{{N}\choose{3}}+ {{N}\choose{5}} + {{N}\choose{7}}}
	 5^{{N}\choose{3}} \\
	 (5,1^2) & \text{\scriptsize{$N(N-1)(N-2)(N+1)(N+2)(N+3)(N+4)/336$}} &
	 3^{{{N}\choose{3}}+{{N}\choose{5}}} 
	 7^{{{N}\choose{5}}} 
	 \\
	 (4,3) & \text{\scriptsize{$N^2(N-1)(N+1)^2(N+2)(N+3)/360$}} & 
	 3^{{{N}\choose{3}} + {{N}\choose{7}}} 5^{{{N}\choose{5}} + {{N}\choose{7}}} \\
	 (4,2,1) & \text{\scriptsize{$N^2(N-1)(N-2)(N+1)(N+2)(N+3)/144$}} & 
	 2^{{{N}\choose{7}}}
	 3^{{{N}\choose{3}}+{{N}\choose{5}}+{{N}\choose{7}}}
	 \\
	 (4,1^3) & \text{\scriptsize{$N(N-1)(N-2)(N-3)(N+1)(N+2)(N+3)/252$}} & 21^{{{N}\choose{5}}} \\
	 (3^2,1) & \text{\scriptsize{$N^2(N-1)(N-2)(N+1)^2(N+2)/240$}} & 
	 2^{{{N}\choose{3}}+{{N}\choose{7}}}
	 3^{{{N}\choose{5}}+{{N}\choose{7}}}
	\\
	 (3,2^2) & \text{\scriptsize{$N^2(N-1)^2(N-2)(N+1)(N+2)/240$}} & 
	 3^{{{N}\choose{5}}+{{N}\choose{7}}}
	 5^{{{N}\choose{7}}}
	\\
	 (3,2,1^2) & \text{\scriptsize{$N^2(N-1)(N-2)(N-3)(N+1)(N+2)/144$}} & 
	 2^{{{N}\choose{7}}}
	 3^{{{N}\choose{5}}+{{N}\choose{7}}}
	\\
	 (3,1^4) & \text{\scriptsize{$N(N-1)(N-2)(N-3)(N-4)(N+1)(N+2)/336$}} & 
	 30^{{{N}\choose{5}}} 105^{{{N}\choose{7}}}
	\\
	 (2^3,1) & \text{\scriptsize{$N^2(N-1)^2(N-2)(N-3)(N+1)/360$}} & 
	 15^{{{N}\choose{7}}}
	\\
	 (2^2,1^3) & \text{\scriptsize{$N^2(N-1)(N-2)(N-3)(N-4)(N+1)/360$}} & 
	 3^{{{N}\choose{7}}}
	\\
	 (2,1^5) & \text{\scriptsize{$N(N-1)(N-2)(N-3)(N-4)(N-5)(N+1)/840$}} & 
	 7^{{{N}\choose{7}}}
	\\
	 (1^7) & {N}\choose{7}  & 35^{{{N}\choose{7}}} \\
\hline
\end{flexitab}
$$

\section{Unitary determinants of the Suzuki group} \label{Suz} 

To illustrate the use of the determinants of the symmetrizations we give 
a small example. 
The covering group $6.Suz$ of the Suzuki group has a 
unitary representation of degree 12 over the field $K=\Q[\sqrt{-3}]$. 
This representation fixes a Hermitian form $B$ of determinant 1. 
The symmetrizations of this 12-dimensional Hermitian form up to
degree 7 
allow to find the determinants of quite a few 
irreducible unitary characters of covers of the Suzuki group. 
Note that ${{12}\choose{k}} $ is even for $k=1,2,3,5,6,7,9,10,11$  and 
odd for $k=4,8$. 

$$
\begin{array}{cccc}
\mbox{ partition } & \chi &  \chi(1) & \mbox{ det } \\
\hline
(1,1) & 78 & 66  & 1 \\
( 2 ) & 80 & 78 & 1 \\
( 2, 1 ) & 48 & 572 & 1 \\
( 3 ) & 46 & 364 & 1 \\
( 2^2 ) & 85 & 1716 & 3 \\ 
( 1^5 ) & 154+156 & 12+780 & 1 \times 1 \\
( 2, 1^3 ) & 156 + 162 & 780 + 4368 & 1 \times 1 \\
( 2^2, 1 ) & 164 & 8580 & 1 \\
( 3, 1^2 ) & 174 & 12012 & 1 \\
(3,2) & 170 & 12012 & 1 \\
( 4, 1 ) & 172 & 12012 & 1 \\
( 5 ) & 160 & 4368 & 1 \\
( 4, 1^2 ) & 21 & 50050 & 1 \\
(1^7) & 153 + 155 & 12 + 780 & 1 \times 1 \\
(2,1^5) & 153+155+157+163 & 12+780 + 924 + 8580 & 1 \\
(2^2,1^3) & 155 + 161 + 163 + 179 & 780+4368+8580 + 27456 & 1 \\ 
(2^3,1) & 155 + 161 + 163 + 185 & 780 + 4368 + 8580 + 42900 & 1 \\
(3,1^4) & 157+175+179 & 924 + 23100 + 27456 & 1 \\
(3,2,1^2) & 163 +179 + 199 & 8580 + 27456 +144144 & 1 \\
(3,2^2) & 173 + 185+191 & 12012 + 42900 + 77220 & 1 \\
(3^2,1) & 169+199 & 12012 + 144144 & 1 \\
(4,1^3) & 175+193 & 23100 + 105600 & 1 \\
(4,2,1) & 207 & 300300 & 1 \\ 
(4,3) & 169+201 & 12012+144144 & 1 \\
(5,1^2) & 203 & 171600 & 1 \\
(5,2) & 171+205 & 12012 + 180180 & 1 \\
(6,1) & 159+195 & 4368 + 112320 & 1 \\
(7) & 159+177 & 4368+27456 & 1 
\end{array}
$$
\newpage

Just considering the faithful characters 
(with numbers 153 to 210) of $6.Suz$ this allows us to 
conclude that all the characters that occur in one of the symmetrizations 
above of one of the two complex conjugate characters $153$ or $154$ 
have unitary determinant 1. 
The missing ones are 
$ 165-168,
181-184, 187-190 $ which have  character fields of 
degree 4 over the rationals, and 
$197, 198$ and $209, 210$. 
The latter four characters occur in $\Sym_{\lambda }(\chi_{153})$ 
resp. $\Sym_{\lambda }(\chi_{154})$ for 
$\lambda = (3,2^4)$ resp. $\lambda = (4,2,1^5)$ which 
are sums of even degree absolutely irreducible characters with
character field $\Q[\sqrt{-3}]$. 

The complex conjugate characters $\chi_{21}$ and $\chi _{22}$ 
are the unique absolutely irreducible characters of the simple
Suzuki group $Suz$ of even degree and
indicator $o$. They have degree $50050$ and the computation above shows that
their determinant is 1. 
Similary we can conclude that all (4 pairs of)
even degree indicator $o$ characters
of $2.Suz$  have determinant 1, where we need to consider 
$$\lambda = (4,1^5) \mbox{ of dim. } 1/6480 \prod_{i=-5}^3 (N+i) 
\mbox{ and det } 3^{{N}\choose{7}} $$
and
$$\lambda = (3^3) \mbox{ of dim. } 1/8640 (N-2)(N-1)^2N^3(N+1)^2(N+2) 
\mbox{ and det } 
6^{{N}\choose{3}} 10^{{N}\choose{7}} 15^{{N}\choose{9}} .
$$

For the faithful characters of $3.Suz$, these symmetrizations  
show the ones of degree $66, 78,$ and $1716$ have determinant 1.

\section{Refined symmetrizations for orthogonal groups} 

In this section we assume that $\Char(K) = 0$ and 
that $B:V\times V\to K$ is a symmetric bilinear form. 
Then there are 
 ${n}\choose{2}$ linearly independent $O(B)$-invariant
 epimorphisms $\pi _{ij}$ for all  $1\leq i < j\leq n $
 by evaluating $B$ in positions $i,j$ of 
 the tensors: 
  $$\pi _{ij} : \otimes ^n V \to \otimes ^{n-2} V, 
  \pi_{ij} (w_1\otimes \ldots \otimes w_n ) = B(w_i,w_j) 
  w_1\otimes \ldots \otimes w_n  $$ 
  where in the last tensor product the vectors $w_i$ and $w_j$ are 
  omitted. 
  There are $O(B)$-invariant monomorphism $ \varphi _{ij} :\otimes ^{n-2} V \to 
  \otimes ^n V  $ with $\pi_{ij} \circ \varphi _{ij} = N \id _{\otimes ^{n-2} V } $: Choose a basis 
  $(v_1,\ldots , v_N)$  of $V$ and define
  $(v_1',\ldots , v_N') $ to denote  its dual basis, i.e. 
  $B(v_i,v_j') = \delta _{ij} $. Then 
  for all $w_1,\ldots , w_{n-2} \in V$ we put 
  $$\varphi_{ij}( w_1\otimes \ldots \otimes w_{n-2}) := 
  \sum_{k=1}^N (w_1 \otimes \ldots \otimes v_k \otimes \ldots \otimes v_k' \otimes  w_{n-2} ) $$
  where $v_k$ is inserted in position $i$ and $v_k'$ in position $j$. 
 The compositions $\varphi_{ij} \circ \pi_{ij} $
 are $O(B)$-invariant endomorphism of $\otimes ^n V$ giving rise to 
 generators of the Brauer algebra, the endomorphism algebra of 
 the $O(B)$-module $\otimes ^n V$.
See also \cite{SchurWeylOrthogonal} for an application of these ideas 
to a Schur-Weyl-duality for orthogonal groups in odd characteristic.

  \begin{proposition} 
	  For all $1\leq i < j \leq n$ the maps 
	  $$\varphi _{ij} :(\otimes ^{n-2} V , N\otimes ^{n-2} B) \to
	  (\otimes ^n V , \otimes^n B)  $$
	  are $O(B)$-invariant isometric embeddings. 
  \end{proposition} 

  \begin{proof}
	 For a pair $(v_1,\ldots , v_N) $ and $(v'_1,\ldots , v'_N )$ of
	  dual basis of $V$ the Gram matrices are inverse to each other,
	  so 
	  $$((B(v_i,v_j)_{1\leq i,j\leq N}) 
	  ((B(v'_i,v'_j)_{1\leq i,j\leq N})  = I_N .$$
	  Now observe that 
	  $$\otimes ^2B (\sum_{k=1}^N v_k\otimes v'_k, \sum_{k=1}^N v_k\otimes v'_k )  = \sum _{k,\ell } B(v_k,v_{\ell }) B(v'_k,v'_{\ell }) 
	   = \sum _{k,\ell } B(v_k,v_{\ell }) B(v'_{\ell },v'_k) =N  $$ 
	   as the trace of the product of these two Gram matrices.
  \end{proof}

  Note that the same proof also works for a non-degenerate 
skew-symmetric bilinear form $B$, 
  where interchanging $v'_k$ and $v'_{\ell}$ introduces a minus sign, 
and hence we need to replace $N$ by $-N$ in the formula of the
proposition. 

\begin{remark}
For the full symmetrizations there was an orthogonal decomposition 
according to the different $n$-element multi-subsets 
	$$C= \{ i_1^{x_1},\ldots , i_k^{x_k}\} 
	\subset \{ 1,\ldots , N \}$$
where the dimension of the orthogonal summand with Gram matrix $X_C$ 
	did only depend on the composition $(x_1,\ldots , x_k)$ 
	of $n$.
For the refined symmetrization, the dimensions of the
	naturally occurring orthogonal summands  grow with 
	$N=\dim (V)$.
	Following 
 \cite{Frame}, who uses \cite{MN} 
	to determine the absolutely irreducible $O(B)$-submodules 
	of $\Sym_{\lambda }(V)$ for partitions $\lambda $ of $n$ 
	and $n\leq 6$, we write
$$\Sym_{\lambda }(V) \cong \Sym_{\lambda }'(V) \oplus \bigoplus _{\gamma}
	m(\lambda , \gamma ) \Sym_{\gamma }'(V) $$
	where $\gamma $ runs over certain partitions of $n-2$, $n-4$, $\ldots$. 
	Here $m(\lambda, \gamma )\in \N$ is the multiplicity of 
	$\Sym_{\gamma }'(V)$ as a composition factor of
	$\Sym_{\lambda }(V)$. 
As $\Sym_{\gamma} '(V)$ are absolutely irreducible, there is a one-dimensional
	space of $O(B)$-invariant quadratic forms on these modules. 
So there are symmetric invertible matrices 
	$$c(\lambda , \gamma ) \in K^{m(\lambda,\gamma)\times m(\lambda , \gamma ) }$$ such that
$$\Sym_{\lambda }(B) \cong \Sym_{\lambda }'(B) \oplus \bigoplus _{\gamma}
	 c(\lambda , \gamma) \otimes \Sym_{\gamma }'(B) .$$
To determine the values of $c(\lambda ,\gamma )$ 
for a partition $\gamma $ of $n-2$ 
(with $m(\lambda , \gamma ) = 1 $) 
it suffices to 
choose a suitable embedding $\varphi_{ij}:\otimes ^{n-2}V \to \otimes ^{n} V $ and compute 
$$c(\lambda , \gamma ) = 
(\otimes ^n B(e_{\lambda } \varphi_{ij}(v), 
e_{\lambda } \varphi _{ij} (v) )  ) / (\otimes ^{n-2} B(v,v) ) $$ 
for a suitable $v\in \Sym_{\gamma }'(V) $, e.g. $v=e_{\gamma } (v_1\otimes 
	\ldots \otimes v_{n-2})$.
For partitions $\gamma $ of smaller $n$ we need to iterate this procedure, 
i.e. consider
$$\varphi_{k,l} \circ \varphi_{i,j} : \otimes ^{n-4}V \to \otimes ^{n} V$$
and so on. 
\end{remark}

  The following sections contain some examples illustrating an 
  elementary way to compute determinants of refined symmetrizations
  for non-degenerate symmetric bilinear forms. 
  As in the previous sections we fix an orthogonal basis 
  $(v_1,\ldots ,v_N)$ of $V$ and put $a_i:=B(v_i,v_i)$.

\subsection{The refined symmetrizations for $(1^n)$ and  $(n)$}

\begin{remark}
For the exterior power all
 $\Sym_{(1^n)} (V) _{ij} $ are zero.
As $\dim (\Sym_{(1^n)}(V)  ) = {{N}\choose{n}} $, this dimension
is even, if and only if $c((1^n),N) =  (n!)^{{N}\choose{n}} $ 
is a square. 
Then 
$$\det(\Sym_{(1^n)}(B)) = \det(B) ^{{{N}\choose{n} } \frac{n}{N}}
= \det(B) ^{{{N-1}\choose{n-1} }} .$$
\end{remark}

\begin{proposition}
$$\Sym_{(n)}(V) \cong  
	\bigoplus_{j=0}^{\lfloor n/2 \rfloor}  
	\Sym_{(n-2j)}'(V) $$ 
	with 
	$$c((n),(n-2)) = 2 n! (n-2)!  (N+2(n-2)) $$ 
	and
	$$c((n),(n-4)) = 8 n! (n-4)! ( N + 2(n - 4) ) ( N + 2(n - 3) ).$$
\end{proposition}

\begin{proof}
Let $I=\{1,\ldots , n-2 \}$ and $w:=e_{n-2} (v_1\otimes \ldots \otimes 
	v_{n-2} ) \in \Sym_{(n-2)}'(V)$. 
	Then $\otimes ^{n-2}B(w,w) = (n-2)! q(I)$ and 
	$$
	\begin{array}{l} 
		e_{(n)} \varphi_{12}(w) = 
	\sum _{k=1}^{N} \frac{1}{a_k} e_{(n)} v_k\otimes v_k \otimes w  = \\
		(n-2)! \sum_{k=1}^N \frac{1}{a_k} e_{(n)} 
		v_k\otimes v_k \otimes v_1 \otimes \ldots \otimes v_{n-2} =\\
		(n-2) !	\sum _{k=1}^{n-2} \frac{1}{a_k} e_{(n)} v_k\otimes v_k \otimes  
	v_1 \otimes \ldots \otimes v_{n-2} \\
		+(n-2) !\sum _{k=n-1}^{N} \frac{1}{a_k} e_{(n)} v_k\otimes v_k \otimes  
	v_1 \otimes \ldots \otimes v_{n-2}
	.
	\end{array}$$ 
The norms of the pure tensors in the first sum are 
	$3! n! q(I)$ whereas the pure tensors in the second sum have norm 
	$2! n! q(I)$. 
So in total we compute 
	$\otimes^nB(e_{(n)}\varphi_{12}(w) ) =$
	$$
	(n-2)!^2((n-2) 3! n! +(N-(n-2)) 2! n!) q(I) = 
	(n-2)!^2 2 n! (N+2(n-2))  q(I)$$
	and hence 
	$$c((n),(n-2)) = 2 (n-2)! n! (N+2(n-2)) = (n-2)!^2 (2(n(n-1))) (N+2n-4).$$
	To compute $c((n),(n-4))$ we choose $I=\{1,\ldots , n-4\}$ 
	and $w:= e_{(n-4)} v_1\otimes \ldots \otimes v_{n-4} \in 
	\Sym_{(n-4)}'(V)$ of norm $(n-4)! $. 
	Now 
	$$e_{(n)}\varphi_{12} \circ \varphi_{34} (w) =
	e_{(n)} \sum_{k=1}^N \sum_{\ell =1}^N \frac{1}{a_k} \frac{1}{a_{\ell}} 
	v_k \otimes v_k \otimes v_{\ell } \otimes v_{\ell } \otimes w $$
	has summands of 5 different types: 
	$$
	\begin{flexitab} 
		{
                        type = array,
                        stretch=1.1,
                        columns=lcc
                        }
			\mbox{ case } & \mbox{ norm / norm(w) } & \mbox{ anz } \\ \hline 
		k=\ell \in I & n! (n-4)! 5! & n-4 \\ 
		k\neq \ell \in I & n! (n-4)! 3!3!2 & (n-4)(n-5) \\ 
		k\in I, \ell \not\in I & n! (n-4)! 3!2!2 & (n-4)(N-(n-4)) \\ 
		k\not \in I, \ell \in I & n! (n-4)! 2!3!2 & (n-4)(N-(n-4)) \\ 
		k= \ell \not\in I & n! (n-4)! 4! & (N-(n-4)) \\ 
		k\neq \ell, k,\ell \not\in I & n! (n-4)! 2!2!2 & (N-(n-4)) (N-(n-3)) 
		\\ 
		\hline
	\end{flexitab}
	$$
	So $$\begin{array}{l} c((n),(n-4)) =  \\
	n! (n-4) ! ( 5! (n-4) + 2\cdot 3!^2(n-4)(n-5) + \\
		2^2(2!3!(n-4)(N-(n-4))) +4! (N-(n-4)) + 8 (N-(n-4))(N-(n-3)) =
		\\
		n! (n-4)! 8( N + 2n - 8 ) ( N + 2n - 6 ) \\
	\end{array} $$
\end{proof}

\begin{remark}
Of course it is possible to continue like this and compute 
	$c((n),(n-2j))$ for $j\geq 3$. 
However this becomes more and more tedious, I think that there 
should be a better way, as 
	$\Sym_{(n)}'(V)$ is just the space of harmonic polynomials
	for the Laplace operator 
	$\sum_{i=1}^N a_i \frac{d^2}{dx_i^2} $ 
	associated with the quadratic 
	form  defined by $B$.
\end{remark}

\subsection{The refined symmetrizations for $(2,1)$ and $(3,1)$}

\begin{proposition}
	$$(\Sym_{(2,1)} (V), \Sym_{(2,1)}(B))  \cong 
	(V,8(N-1)B) \perp (\Sym_{(2,1)}'(V), \Sym_{(2,1)}'(B) ) .$$ 
	Hence 
	$\dim (\Sym_{(2,1)}'(V)) = \frac{1}{3}N(N-2)(N+2) $ and 
	$$ \det (\Sym_{(2,1)}'(B)) = 3^{{N}\choose{3}} 2^N (N-1)^N 
	\det(B) ^{N^2-2} .$$
	In particular the dimension of the refined 
	symmetrization $\Sym_{(2,1)}'(V)$ 
	is even, if and only if, $\dim(V) = N$ is even. 
	In this case ${N}\choose{3}$ is also even and hence 
	$\det(\Sym_{(2,1)}'(B)) $ is a square.
\end{proposition}

\begin{proof}
	A basis of 
	$e_{(2,1)} \varphi_{12}(V) $ is given by 
	$$ (b_i:=\sum_{k=1}^N \frac{1}{a_k} e_{(2,1)}  (v_k\otimes v_k \otimes v_i ) 
	\mid 1\leq i\leq N ) .$$
	Now $b_i = \sum_{k\neq i} \frac{2}{a_k} (v_k\otimes v_k \otimes v_i - 
	v_i \otimes v_k \otimes v_k )  $ satisfies 
	$$\otimes ^3B(b_i,b_i) = 8 \sum_{k\neq i} \frac{1}{a_k^2} B(v_k,v_k) B(v_k.v_k) 
	B(v_i,v_i) = 8(N-1) B(v_i,v_i) .$$
\end{proof}

\begin{proposition}
        $$\Sym_{(3,1)} (V)  \cong
        \Sym_{(1,1)}(V) \oplus 
        \Sym_{(2)}'(V) \oplus \Sym_{(3,1)}'(V) ,$$
        so $$\dim (\Sym_{(3,1)}'(V)) = \frac{1}{8} 
        (N-2)(N-1)(N+1)(N+4). $$
        Up to squares we get
        $c((3,1),(1,1)) = 2(N+2)$ and
        $c((3,1),(2)) = N$.
\end{proposition}

\begin{proof}
Put
        $$  b_1 := \sum_{k=1}^N \frac{1}{a_k} e_{(3,1)}  (v_k\otimes v_k \otimes v_1 \otimes v_2)  \mbox{ and } 
          b_2 := \sum_{k=1}^N \frac{1}{a_k} e_{(3,1)}  (v_k\otimes v_k \otimes v_2 \otimes v_1)  .$$
Then
$$\otimes ^4B(b_1,b_1) = 
\otimes ^4B(b_2,b_2) = (24(N-2)+6^2 \cdot 2  + 2^2 \cdot 2 ) a_1a_2
= 8(3N+4) a_1a_2 $$
and
$$\otimes ^4 B(b_1,b_2) = (-8(N-2) -2\cdot 8 \cdot 3) a_1a_2 =
-8(N+4) a_1a_2 .$$
So
 $b_1-b_2 = e_{(3,1)} \varphi_{1,2} (e_{(1,1)} (v_1\otimes v_2))  $
has norm $64 (N+2) a_1a_2$ giving $c((3,1),(1,1)) = 2(N+2)$
and $b_1+b_2 = e_{(3,1)} \varphi_{1,2} (e_{(2)} (v_1\otimes v_2))  $
has norm $32 N a_1a_2$ which yields
$c((3,1),(2)) = N$.
\end{proof}

\subsection{The refined symmetrizations for $(2,1^{n-2})$}

\begin{proposition}
	$\Sym_{(2,1^{n-2})} (V) \cong 
	 \Sym_{(1^{n-2})}(V) \oplus \Sym_{(2,1^{n-2})}'(V)$, 
	 so 
	 $$\dim (\Sym_{(2,1^{n-2})}'(V) = 
	 (n-2) {{N}\choose{n}} +
	 (n-1) {{N}\choose{n-1}}.$$ 
	 Up to squares 
	 $$c((2,1^{n-2}),(1^{n-2})) =(n-1) (N-(n-2)).$$
\end{proposition}

\begin{proof}
	Let $I=\{1,\ldots , n-2\}$ and put 
	 $$w := e_{(1^{n-2})} v_{1} \otimes \ldots \otimes v_{n-2} \in 
	 \Sym_{(1^{n-2})}(V) .$$
	 Then the norm of $w$ is $(n-2)!$.
	 For $\lambda := (2,1^{n-2})$ we compute 
	 $$ 
	 \begin{array}{l} e_{\lambda } \varphi_{1,2}(w) = 
	\sum_{k=1}^N \frac{1}{a_k} e_{\lambda}  
	v_k\otimes v_k \otimes w   =  \\
	(n-2)! \sum_{k=1}^N \frac{1}{a_k} e_{\lambda} v_k \otimes v_k \otimes 
	v_1 \otimes \ldots \otimes v_{n-2} = \\
		 (n-2)! \sum_{k=1}^{n-2} \frac{1}{a_k} e_{\lambda} v_k \otimes v_k \otimes 
	v_1 \otimes \ldots \otimes v_{n-2}  + \\
	(n-2)! \sum_{k=n-1}^N \frac{1}{a_k} e_{\lambda} v_k \otimes v_k \otimes 
	v_1 \otimes \ldots \otimes v_{n-2} .
	 \end{array}
		 $$
		 Now the first summand is 0 whereas the last $(N-(n-2))$
		 summands have norm 
		 $(n-2)!^2 2^2 (n-1)! q(I) $. 
		 Hence 
		 $c(\lambda, (1^{n-2}) ) = 4 (n-1)! (n-2) ! (N-(n-2))=(n-1) (2(n-2))^2 (N-(n-2))$. 
\end{proof}

\subsection{Determinants of refined symmetrizations} \label{AUSorth}

This section gives  tables of the $c(\lambda , \gamma )$ 
for partitions $\lambda $ of $n$ and all $n\leq 6$.
Here $(V,B)$ is a non-degenerate symmetric bilinear
space over a field of characteristic 0, $N=\dim (V)$ is assumed to
be $\geq n$.
For fixed partition $\lambda $ of $n$ 
we display the partitions $\gamma $ of $m\leq n$ 
and  give $c(\lambda ,\gamma )\in \N$ (up to squares) such that 
$$\Sym_{\lambda }(B) \cong \Sym_{\lambda }'(B) \oplus _{\gamma } 
c(\lambda , \gamma ) \Sym_{\gamma }'(B). $$
We omit the rows for $1^n$ since $\Sym_{(1^n)}(V) = \Sym_{(1^n)}'(V)$.
The composition factors of $\Sym_{\lambda }(V)$ are taken from the 
table in \cite[p. 157]{Frame}.
The only composition factor that occurs with multiplicity $>1$ is
$\Sym_{(2)}'(V)$ in $\Sym_{(4,2)}(V)$, where the 
multiplicity is 2.
Here 
$$A:=c((4,2),(2)) = 2^8 \left( \begin{array}{rr} 
	12N^2 - 32    &     2N^2 - 32 \\
	      2N^2 - 32 & 7N^2 + 20N - 32
\end{array} \right) 
	\in \Z[N]^{2\times 2} $$ 
is of determinant 
$$ 2^{20} 5 (N-2)N(N+1)(N+4) .$$
To compute $A$ we chose the embeddings 
$$f_1:= e_{(4,2)} \circ \varphi _{1,2} \circ \varphi _{3,4} \mbox{ and } 
f_2 := e_{(4,2)} \circ \varphi _{1,2} \circ \varphi _{5,6} : \Sym_{(2)}'(V) \to \Sym_{(4,2)}(V) .$$ 
For $v=e_{(2)}(v_1\otimes v_2) \in \Sym_{(2)}'(V)$ 
we computed $A_{ij} = \otimes^6B(f_i(v) , f_j(v)) .$

$$
\begin{flexitab}{
                type = array,
                stretch = 1.1,
                columns = |c|c|l|
                }
		\hline
\lambda  &  \gamma  & c(\lambda ,\gamma )  \\ 
\hline
	(2) & () & N \\
\hline
	(3) & (1) & 3(N+2) \\
 (2,1) & (1) & 2(N-1) \\
\hline
	(4) & (2) ,\ () & 6(N+4) ,\ 3N(N+2) \\
	(3,1) & (2),\ (1^2) & N,\ 2(N+2) \\
	(2^2) & (2),\ ()  & N-2,\ 2N(N-1) \\
  (2,1^2) & (1^2) & 3(N-2) \\
\hline
	(5) & (3) ,\ (1) & 10(N+6),\ 15(N+2)(N+4) \\
	(4,1) & (3),\ (2,1),\ (1) & 6(N+1) ,\ 5(N+4) ,\ 6(N-1)(N+2) \\
	(3,2) & (3),\ (2,1),\ (1)  & (N-2),\ 3(N+1),\ (N-1)(N+2) \\
	(3,1^2) & (2,1),\ (1^3) & 6(N-1),\ 15(N+2) \\
	(2^2,1) & (2,1),\ (1) & 6(N-3),\ 6(N-1)(N-2) \\
 (2,1^3) & (1^3) & (N-3) \\
\hline
	(6) & (4),\ (2) & 15(N+8) ,\ 5(N+4)(N+6)  \\
	 & () & 15N(N+2)(N+4) \\
	(5,1) & (4),\ (3,1) &  2(N+2),\ N+6 \\
	   & (2),\ (1,1) &   6N(N+4) ,\ (N+2)(N+4) \\
	(4,2) & (4),\ (3,1),\ (2,2) & N -2,\ 10(N+2) ,\ 30(N+4) \\
	& 2\cdot (2) ,\ ()  & A ,\ N(N-1)(N+2)\\
	(4,1,1) &  (3,1),\ (2,1,1),\ (1,1)& N ,\ 2(N+4),\ (N-2)(N+2) \\
	(3,3) & (3,1),\ (1,1) & 3N ,\ 3(N+1)(N+2) \\
	(3,2,1) & (3,1),\ (2,2),\ (2,1,1) & 6(N-3),\ 2(N-1),\ 15(N+1)  \\
	& (2),\ (1,1) &  N(N-2) , 15(N+2)(N-2) \\
	(2^3) & (2,2) ,\ (2),\ () & N-4,\ 2(N-2)(N-3) , \ 6N(N-1)(N-2) \\
	(3,1^3) & (2,1^2),\ (1^4) & 2(N-2) ,\ 6(N+2) \\ 
	(2^2,1^2) & (2,1^2),\  (1^2) & 2(N-4), \ 3(N-2)(N-3) \\ 
	(2,1^4) & (1^4) & 5(N-4)  \\
	\hline
\end{flexitab}
$$

Note that the same programs prove that the refined symmetrization 
$\Sym_{\gamma }'(V) $ does occur as an orthogonal summand of the 
$O(B)$-module $\Sym_{\lambda }(V)$, whenever we compute $c(\lambda , \gamma ) 
\neq 0$. We can also obtain the multiplicity (at least a lower bound $m$) 
by finding $c(\lambda, \gamma ) \in K^{m\times m}$ of full rank $m$. 

\subsection{An example: the refined symmetrizations of 
the 24-dimension representation of $2.Co_1$} 

The covering group of the sporadic simple Conway group $2.Co_1$ 
is a subgroup of index 2 of the automorphism group of 
the 24-dimensional extremal unimodular lattice, the Leech lattice. 
So this group has a 24-dimensional rational representation of 
determinant 1. 
The refined symmetrizations for partitions of even numbers 
hence yield representations of the simple group $Co_1$. The
table below lists the ones that are orthogonally stable together
with their decomposition into irreducibles and the corresponding 
determinants as predicted by the table in Section \ref{AUSorth}. 
For partitions of odd numbers, the refined symmetrizations are
faithful representations of $2.Co_1$ and hence have 
orthogonal determinant 1 (see for instance \cite[Theorem 4]{survey})
which is confirmed by the computation of the determinants of 
the refined symmetrizations. 
\\

\begin{center} 
{\bf Determinants of orthogonally stable characters for $Co_1$:}  \\
$
\begin{array}{cccc}
\mbox{ partition } & \chi &  \chi(1) & \mbox{ det } \\
\hline
	(2) & 2 & 276 & 1 \\
	(2,1,1) & 8 & 37674 & 1 \\
	(2,2) & 7 & 27300 & 253 \\
	(4) & 6 & 17250 & 91 \\ 
	(3,1^3) & 15+23 & 483000 + 1771000 & 1\\
	(3,2,1) & 31 & 4100096 & 161 \\
	(5,1) & 24 & 1821600 & 65 \\
	(4,2) & 29 & 2816856 & 13 \\
	(2^3) & 12+19 & 313950+822250 & 77 \\
	(6) & 10+14 & 80730 + 376740 & 35 
\end{array}
$
\end{center}

\end{document}